\documentclass[11pt]{article}
\usepackage{doublespace}
\usepackage{amsmath}
\usepackage{amssymb}
\usepackage[all]{xy}
\newcommand{\hide}[1]{#1}
\newcommand{\bea}[1]{\begin{eqnarray}\label{#1}}
\newcommand{\eea}{\end{eqnarray}}
\newcommand{\beas}{\begin{eqnarray*}}
\newcommand{\eeas}{\end{eqnarray*}}
\newcounter{rocount}
\newenvironment{rolist}{\begin{list}{(\roman{rocount})}{\usecounter{rocount}}}{\end{list}}
\newenvironment{Rolist}{\begin{list}{(\Roman{rocount})}{\usecounter{rocount}}}{\end{list}}

\newcounter{gracount}

\newcommand{\agnorm}[2]{\norm{#1}_{#2}}

\newtheorem{theor+}{Theorem}[section]
\newtheorem{propo+}[theor+]{Proposition}
\newtheorem{corol+}[theor+]{Corollary}
\newtheorem{lemma+}[theor+]{Lemma}
\newtheorem{obser+}[theor+]{Observation}
\newtheorem{defin+}[theor+]{Definition}
\newtheorem{examp+}[theor+]{Example}
\newtheorem{remar+}[theor+]{Remark}
\newtheorem{conje+}[theor+]{Conjecture}
\newtheorem{quest+}[theor+]{Question}
\newtheorem{pictu+}[theor+]{Picture}
\newtheorem{notes+}[theor+]{Notes}
\newtheorem{exerc+}[theor+]{Exercise}
\newtheorem{backg+}[theor+]{Background}
\newtheorem{fact+}[theor+]{Fact}
\newenvironment{theor}[1]{\begin{theor+}\label{#1}\slshape}{\end{theor+}}
\newenvironment{propo}[1]{\begin{propo+}\label{#1}\slshape}{\end{propo+}}
\newenvironment{corol}[1]{\begin{corol+}\label{#1}\slshape}{\end{corol+}}
\newenvironment{lemma}[1]{\begin{lemma+}\label{#1}\slshape}{\end{lemma+}}

\newenvironment{defin}[1]{\begin{defin+}\label{#1}\upshape}{\end{defin+}}

\newenvironment{remar}[1]{\begin{remar+}\label{#1}\upshape}{\end{remar+}}
\newenvironment{conje}[1]{\begin{conje+}\label{#1}\upshape}{\end{conje+}}

\newenvironment{demo}{{\bf\sl Proof:}}{\hfill$\Box$}

\renewcommand{\Box}{\blacksquare}
\newcommand{\NN}{{\mathbb{N}}}
\newcommand{\ZZ}{{\mathbb{Z}}}

\renewcommand{\epsilon}{\varepsilon}
\renewcommand{\phi}{\varphi}
\newcommand{\norm}[1]{{\left\|{#1}\right\|}}

\newcommand{\cl}{^{\overline{\makebox[0.4em]{\rule{0em}{0.6ex}}}}}

\newcommand{\infer}{\mbox{$\Longrightarrow$}}

\newcommand{\cstar}{\mbox{$C^*$}}
\newcommand{\sa}{_{sa}}
\newcommand{\dd}{^{**}}
\newcommand{\po}{_+}

\newcommand{\id}{\operatorname{id}}

\newcommand{\cco}[1]{\overline{\operatorname{co}}({#1})}

\newcommand{\gkp}[1]{\cite[{#1}]{gkp:cag}}
\newcommand{\CCC}{{\mathbb K}}
\newcommand{\BBB}{{\mathbb B}}

\newcommand{\myunit}{\text{\sf 1}}
\newcommand{\myw}{z}

\newcommand{\wea}{\mbox{\sl weak}}
\newcommand{\weasta}{\mbox{\sl weak$^*$}}
\newlength{\iteid}

\newcommand{\cantor}{M}

\newcommand{\HHH}{{H}}

\newcommand{\MMM}{{\mathfrak M}}

\newcommand{\newsym}[1]{}

\newcommand{\smalltwobytwo}[4]{\left[\begin{smallmatrix}{#1}&{#2}\\{#3}&{#4}\end{smallmatrix}\right]}

\newcommand{\inner}[2]{\left\langle{{#1},{#2}}\right\rangle}

\newcommand{\matrM}{\mathbf M}

\newcommand{\cref}[1]{Conjecture {#1}}

\newcommand{\myz}{\pmb{\mathsf z}}
\newcommand{\PPA}[1]{P({#1})}           %
\newcommand{\A}{{\mathfrak A}}
\newcommand{\I}{{\mathfrak I}}
\newcommand{\J}{{\mathfrak J}}
\newcommand{\PPP}[1]{P({#1})}           %
\newcommand{\pbar}{{\overline{p}}}
\newcommand{\myu}{{1}}
\newcommand{\KKK}{\mathbb{K}}
\newcommand{\SSS}{S}%
\newcommand{\sss}{s}%
\newcommand{\TTT}{\Omega}%
\newcommand{\ttt}{\omega}%
\newcommand{\wlim}{\text{\textsl{weak}}\lim}

\title{Regularity of projections revisited} %
\date{\today}
\author{Charles A.\ Akemann \and  S\o ren  Eilers}

\begin{document}
\maketitle
\begin{abstract}
  The concept of \emph{regularity} in the meta-topological setting of
  projections in the double dual of a \cstar-algebra addresses the
  interrelations of a projection $p$ with its closure
  $\overline{p}$, for instance in the form that such projections act
   identically, in norm, on elements of the \cstar-algebra. 
This concept has been given new actuality with the recent plan of Peligrad
and 
Zsido to find a meaningful notion of Murray-von Neumann type
equivalence among open projections. 

Although automatic in the commutative case, it has been known 
since the late sixties that regularity fails for many
  projections. The original investigations, however, did not answer a
  question such as:
\emph{Are all open and dense projections regular in $\A$, when $\A$ is simple?}
We report here that this and related questions have negative answers.
In the other direction,
we supply 
positive results on regularity of large open projections.
\end{abstract}

%
%
%

\section{Introduction}
For a pair of elements $x,b$ in a \cstar-algebra $\A$ with $b\geq 0$ and $\norm{b}=1$, consider the
quantity 
\[
\agnorm{x}{b}=\sup_{n\geq0} \norm{xb^{1/n} x^*}^{1/2}
\]
It is easy to see (Lemma \ref{agisnorm} below) that $\agnorm{\cdot}{b}$ is a
seminorm, and that $\agnorm{\cdot}{b}\leq \norm{\cdot}$, but
clearly $\agnorm{\cdot}{b}$ will fail to be a norm if $b$ is a zero-divisor.
In fact, for many \cstar-algebras $\A$, $\agnorm{\cdot}{b}$ coincides 
with $\norm{\cdot}$ for all positive elements $b$ of $\A$ that have norm one and are not
zero-divisors. This is easily seen  when $\A$ is Abelian,
and we observe below (Proposition \ref{agmatrix}) that it also is the case when $\A$
is of the form $C(X,\matrM_n)$.

Prompted by a question of Peligrad and Zsido (\cite{cplz:opccr}) we 
answer here, in the negative, the question of whether
$\agnorm{\cdot}{b}$ is a norm for a positive, norm one non-zero
divisor $b$ in a \cstar-algebra $\A$.
Intriguingly, at the present stage of our investigations,
we can not provide a direct, self-contained example of a positive,
norm one
element $b$ in a \cstar-algebra which is not a zero-divisor, and for
which $\agnorm{\cdot}{b}\not= \norm{\cdot}$. Nevertheless, the
theory of \emph{open dense nonregular projections}, to which we
contribute in this paper, shows that such elements do exist in, e.g.,
the stabilized $2^\infty$ UHF algebra.

The concept of \emph{regularity} for a projection in the double dual of a
\cstar-algebra $\A\dd$ was introduced by Tomita in \cite{mt:stoai}
and subsequently studied by Effros and Akemann in 
\cite{ege:oicid},\cite{caa:gswp},\cite{caa:aumac}. 
We are going to collect several known 
characterizations of this notion in Section \ref{regudef} below. If we
restrict our attention to the classes of \emph{open} and \emph{dense}
projections in $\A\dd$ (see Section \ref{notation}), regularity
of the dense and open projection $p$ reduces to the simple property that
\[
\norm{a}=\norm{ap}.
\]
for all $a\in\A$. The relevance of this concept to our discussion of
$\agnorm{\cdot}{b}$ is perhaps already clear to the reader, but will
be explained in depth in Section \ref{onag}.

Every projection in the double dual of a separable \cstar-algebra  dominates
enough regular projections to be a strong limit of a sequence  of such
\cite[Corollary 6]{mt:stoai}. It it clear that every closed projection is
regular, and not hard to prove (a result attributed to
Kaplansky) that every
central projection is regular.
An example of
an open and nonregular projection was found in \cite[Example I.2]{caa:lisc},
and it was proved in \cite{caa:gswp} that every open projection of a von
Neumann algebra is regular, but the initial work of the first author left open the following questions for a
general \cstar-algebra $\A$:
\begin{Rolist}
\item Are all open and dense projections regular in $\A$?
  \item Are all open projections regular in $\A$, when $\A$ is simple?
\item Are all open and dense projections regular in $\A$, when $\A$ is  simple?
\end{Rolist}

We must report here that
there exist $AF$ algebras in which (I)--(III) fail; in fact, problem
(I) is equivalent to the problem of finding an element $b$ as
described above. The existence of such projections is expected to
cause complications 
in the program of \cite{cplz:opccr}  to find a meaningful notion of Murray-von Neumann
type
equivalence among open projections, which sparked a renewed interest
in regularity. However,  our investigations turned up a collection of
positive results, proving that 
many dense and open projections are in fact regular.  It is our hope
that they can be used to work around such  complications.
\section{Preliminaries}
\subsection{Notation}\label{notation}

We are going to work a great deal with the \cstar-algebra $\KKK$ of
compact operators on a separable Hilbert space, and with the $2^\infty$ 
UHF algebra which shall be denoted by
$\MMM$.
We denote generic \cstar-algebras by  $\A$, and always consider them
to be subalgebras of
their double dual $\A\dd$. We denote by $\myz$ the sum of all minimal
central projections of $\A\dd$. This is exactly the central cover of
the  atomic representation of $\A$, i.e. the cover of the sum of all
irreducible representations.

In $\A^*$ we will be working with the sets $Q(\A)$, $S(\A)$,
$P(\A)$ consisting of quasi-states, states and pure states, respectively.
Given a projection $p$, we also need:
\begin{gather*}
L(p)=^\perp\!(\A\dd(1-p))=\{\phi\in\A^*|\forall a\in\A\dd:\phi(a(1-p))=0\}\\
F(p)=\{\phi\in Q(\A)|\phi(1-p)=0\}\\
\PPP{p}=F(p)\cap P(\A)=\{p\in P(\A)|\phi(p)=1\}
\end{gather*}
and call these sets the left invariant subspace, 
the face and the set of pure states {\em supported} by $p$,
respectively.
When $p\in \A\dd$ is a projection, if one of the sets $L(p), F(p),
C(p)$ is \weasta\ closed, so are the other two. This property defines
a \emph{closed} projection (\cite[3.11.9]{gkp:cag}). The \emph{closure} of a
general projection $p$ is the smallest closed projection
$\overline{p}$ which dominates $p$; it is uniquely determined by
\[
L(\overline{p})=\overline{L(p)}.
\]
As we shall see below, it is not true in general that
$F(\overline{p})=\overline{F(p)}$, leading to the definition of
regular projections.

With definitions of closedness and closure in place, we can derive the
following notions in a straightforward way: An \emph{open} projection
is one whose complement is 
closed, and a \emph{dense} projection is one whose closure is the identity.
We also define \emph{compact} projections as those closed projections
which are dominated by an element in $\A$.

\subsection{Regularity}\label{regudef}
Let us collect several known characterizations of regularity as follows:
\begin{theor}{regumain}
Let $p\in \A\dd$ be a projection. Then the following conditions are
equivalent:
\begin{rolist}
\item $\forall a\in \A:\norm{a\pbar}=\norm{ap}.$
\item $\forall a\in M(\A): p\leq a\infer \pbar\leq a$
\item $(L(p)_1)\cl=L(\pbar)_1.$
\item $F(p)\cl=F(\pbar).$
\end{rolist}
\end{theor}
The conditions (iii)--(iv) were used in the original work of Tomita and
Effros on regularity, and their equivalence in presence
of unitality was proved in \cite[6.1]{ege:oicid}. The proof works in
the general case with minor modifications, see the remark below.
The condition (i), and its equivalence with (iii), is found in
\cite[II.12]{caa:gswp}. Condition (ii) is from \cite[Theorem 18]{gkp:scccnt}.
\begin{remar}{nonunitalcomplications}
If
$\phi_\lambda\rightarrow\phi$
\weasta\ in $Q(\A)$, and $\A$ is unital, then
$\norm{\phi_\lambda}=\phi_\lambda(1)\rightarrow\phi(1)=\norm{\phi}$.
This is not true in the non-unital case; for instance a net of states
can converge to \weasta\ to zero. But if $\norm{\phi}=1$, then in all
cases $\norm{\phi_\lambda}\rightarrow \norm{\phi}$.
One needs to use this observation to amend Effros' proof of
\cite[6.1]{ege:oicid} to the non-unital case. But it will not carry
through in establishing equivalence between (i)--(iv) and the property 
\begin{rolist}
  \addtocounter{rocount}{4}
\item $\overline{\{\phi\in (\A^*)\po|\phi(1-p)=0\}}=(\overline{L(p)})_+$
\end{rolist}
which is also considered in that result. In fact, it is shown in
\cite{lgb:nrcpoa} that this condition is strictly weaker than
regularity in the non-unital case. The reader is referred to
\cite{se:cetc} for details about this complication.
\end{remar}
\subsection{Preamble on non-zero divisors}\label{onag}

\begin{defin}{formaldefs}
Let a \cstar-algebra $\A$ be given, and fix $b\in\A$ with $b\geq 0$ and
$\norm{b}=1$. We define
\begin{rolist}
\item $\agnorm{x}{b}=\sup_{n\geq0} \norm{xb^{1/n} x^*}^{1/2}$
\item $\tau(b)=\inf_{x\not=0}{\agnorm{x}{b}}/{\norm{x}}$
\end{rolist}
\end{defin}

When $b\in\A$ is positive and has norm one, we write $[b]$ for its range
projection in $\A\dd$.  We then get:
\begin{lemma}{agisnorm}
We have  $\agnorm{x}{b}=\norm{x[b]}$, so $\agnorm{\cdot}{b}$ is a seminorm. The equality
\[
\norm{\cdot}=\agnorm{\cdot}{b}
\]
holds precisely when $[b]$ is a dense regular projection.
\end{lemma}
\begin{demo}
The norm equality holds as $b^{1/n}\nearrow [b]$ in $\A\dd$, and the second claim is immediate from Theorem \ref{regumain}.
\end{demo}

Clearly  $b$  is a zero
divisor precisely when  $\tau(b) = 0$, and if $b$  is strictly
positive, then  $[b]=1$ and
$\tau(b) = 1$.
In
the case of a commutative $\A$, one easily sees that 0  and  1  are
the only possible values of  $\tau(b)$. This generalizes further:

\begin{propo}{agmatrix}
Let $\TTT$ be a locally compact Hausdorff space. If $b\in C_0(\TTT,\matrM_n)$
is positive, has norm one, and is not a zero-divisor, then $\tau(b)=1$.
\end{propo}
\begin{demo}
Suppose that $\norm{x[b]} < 1$ and $\norm{x}= 1$.  Let  $U$ be an open set
such that for  $\ttt$  in  $U$, $\norm{x(\ttt)} > \norm{x[b]}$.
Let then $\ttt_0$ be chosen so that the number $k$ of nonzero eigenvalues of
$b(\ttt)$ is maximal at $\ttt = \ttt_0$  for  $\ttt  \in  U$. We are going to find 
$r>0$ and a neighborhood $V\subseteq U$ of $\ttt_0$ such that
\[
\operatorname{sp}(b(\ttt))\cap (0,r/2)=\emptyset \text{ for all
  }\ttt\in V.
\]

If $k = 0$,
we achieve this by $r=1$ and $V = U$. When $k>0$, label the distinct nonzero eigenvalues of
$b(\ttt_0)$  as
$s_1 < s_2 < ... < s_k$.  Let  $r$  be the minimum distance between any two
eigenvalues of  $b(\ttt_0)$ (including the 0 eigenvalue).
Using functional calculus with spike functions supported around each
nonzero eigenvalue of $b(\ttt_0)$, we find  for each  $j$ an open neighborhood  $V_j$  of $\ttt_0$
such that for any  $\ttt  \in  V_j$, $b(\ttt)$ has an eigenvalue in the interval
$(s_j - r/4, s_j + r/4)$.
With $V=U\cap V_1\cap\dots\cap V_k$ we have then that
 for  each $\ttt  \in  V$, $b(\ttt)$  has an eigenvalue in each
$(s_j - r/4, s_j + r/4)$. By maximality of  $k$, $b(\ttt)$ can't have an
eigenvalue in the interval  $(0, r/2)$.

        Let $f$ be a continuous function on $[0,1]$ that is $0$ at $0$ and
equal to $1$ on $[r/2,1]$.  Then 
\[
(1-f(b))(\ttt) b(\ttt) = 0
\]
for each $\ttt$ in $V$.
Since $\norm{x(\ttt)[b](\ttt)} < 1$ for each $\ttt$ in $V$, therefore $[b](\ttt)$
is not equal to $1$ for $\ttt$ in $V$, hence $(1-f(b))(\ttt)$ is not zero for such
$\ttt$. Now find a Urysohn function $g$ on $\TTT$ that is $1$ at $\ttt_0$ and $0$
outside $V$.  The operator $\ttt \mapsto g(\ttt) (1-f(b))(\ttt)$ is thus a nonzero
element of $C_0(\TTT,\matrM_n)$ that is orthogonal to $b$.
\end{demo}

Our Theorems \ref{denseopen} and \ref{denseopensimple} below will
show that there exists a
positive,  norm one, element of 
  $C([0,1], \CCC)$ or $\MMM\otimes\KKK$ for which $0<\tau(b)<1$.

\begin{remar}{morework}
For  $\alpha\in [0,1]$, define
\[
\mathcal{R}_t  = \{b  \in  \A^+ \mid \norm{b} = 1, \tau(b) = t\}.
\]
Define 
\[
\mathcal{N} = \{b \in  A^+ \mid \norm{b} = 1, b \text{  is not a zero divisor}\}
\]

Our results above open the discussion of when $\mathcal{N}=\mathcal{R}_1$. 
Related questions, which we shall not attempt to answer here, are
what properties $\mathcal{N}_t$ have (other than automorphism invariance),
whether all $\mathcal{R}_t$ are nonempty when $\mathcal{R}_1\not
=\mathcal{N}$, and whether,
when $b\in \mathcal{R}_t$  for  $t < 1$, we can find  $x  \in
\A_+\backslash \{0\}$  such that
$\sup_n||b^{1/n}x|| = t||x||$.
  \end{remar}

\section{Open and nonregular projections}

The present section contains all our examples establishing existence
of open and nonregular projections under additional assumptions. The
first technical results, collected in Section \ref{vvf} below,
provide the foundation for all the results in Sections \ref{dn} and \ref{dns}.
Section \ref{ninM} does not depend upon \ref{vvf}.

\subsection{A vector-valued function}\label{vvf}
The following notation will be used in all of this section.
 $H$  is an infinite dimensional, separable  Hilbert space. In $H$, we choose a
distinguished unit vector $y$. The space   $\TTT$  is 
some  perfect 
and compact
metric space.  Recall that $\TTT$ is separable and let 
$\SSS$ be a countably infinite dense set in  $\TTT$, enumerated
without repetitions as $\SSS=\{\sss_n\}$.

\begin{lemma}{zerowlim}
  For any $\sss\in \SSS$  there is a norm continuous function
\[
\myw:\TTT\backslash\{\sss\}\longrightarrow H
\]
with the properties
\begin{gather*}
\norm{\myw(\ttt)}=1, \text{ for all }\ttt\in \TTT\backslash \{\sss\}\qquad\\
\end{gather*}
and
\begin{gather*}
  \wlim_{\ttt\longrightarrow\sss}\myw(\ttt)=0.
\end{gather*}
\end{lemma}
\begin{demo}
Choose a basis $\{e_n\}$ for $H$.
Recall that $\TTT$ is metric; we may assume that
$\operatorname{diam}(\TTT)<1$. We consider pairwise overlapping annuli
centered at $\sss$ defined as
\[
A_n=\left\{\ttt\in\TTT\left|
    \frac1{n+2}<d(\ttt,\sss)<\frac1{n}\right\}\right. ;
\]
clearly this covers the paracompact space
$\TTT\backslash\{\sss\}$. Choosing  %
a locally finite partition of unity
for $\TTT\backslash\{\sss\}$ subordinate to the $A_n$ we get
a family $\{\psi_n\}$ with the properties
\begin{rolist}
  \item for each $\ttt\in \TTT\backslash\{\sss\}$, $\psi_n(\ttt)\not=0$ for
    no more than two values of $n\in\NN$.
\item for each $\ttt\in \TTT\backslash\{\sss\}$,
  $\sum_{n=1}^\infty\psi_n(\ttt)=1$.  
\item for each $N$, there exists a neighborhood $U$ of $\sss$ with the
  property that
$\psi_1|U=\cdots=\psi_N|U=0$
\end{rolist}
Let
\[
\myw(\ttt)=\sum_{n\in \NN}\sqrt{\psi_n(\ttt)}e_n.
\]
This is locally a finite sum yielding unit vectors by (i) and (iii). Let $v\in
H$ and
$\epsilon>0$ be given. Choose $N$ such that $\inner{v}{e_n}\leq
\epsilon/2$ when $n\geq N$ 
and $U$ a neighborhood of $\sss$ as in (ii). We then have for any $\ttt\in U$,
by (i) again, that
\[
|\inner{\myw(\ttt)}{v}|=
\left|\sum_{n={N+1}}^\infty\sqrt{\psi_n(\ttt)}\inner{e_n}{v}\right|\leq
|\inner{e_{n_1}}{v}|+
|\inner{e_{n_2}}{v}|\leq 
\epsilon.
\]
Here $n_1,n_2$ are the two values of condition (i).
\end{demo}

\begin{lemma}{imprConstr} Let $H, y, \TTT$ and $\SSS$ be as above. For any
  $\delta > 0$ 
there are  functions 
\[
x:\TTT\longrightarrow H
\qquad
\mu:\SSS\longrightarrow (0,1)
\]
with the properties
\begin{rolist}
\item $\norm{x(\ttt)}=1$ for all $\ttt\in\TTT$.
\item $\norm{y-x(\ttt)}<\delta$ for all $\ttt\in\TTT$.
\item $x$ is norm continuous at any $\ttt\in\TTT\backslash\SSS$.
\item for any  $\sss\in\SSS$,
\[
\wlim_{\ttt \longrightarrow \sss} x(\ttt) = \mu(\sss)x(\sss).
\]
\end{rolist}
\end{lemma}
\begin{demo}
Fix $\eta>0$ such that
\[
{1-\frac1{\sqrt{1+\eta}}+\sqrt{\eta}}<\delta.
\]
Choose a basis for $H$ indexed over $\{0\}\cup\NN\times\NN$, say
$\{e_0,e_{ij}\}$, with $e_0=y$. Set
$H_n=\overline{\operatorname{span}\{e_{ni}\mid i\in\NN\}}$. Choose 
functions $\myw_n:\TTT\backslash \{\sss_n\}\longrightarrow H_n$ according to
Lemma
\ref{zerowlim}, and let 
\[
z_n(\ttt)=\begin{cases}
  \sqrt{\eta 2^{-n}}\myw_n(\ttt)&\ttt\not=\sss_n\\0&\ttt=\sss_n\end{cases} 
\]
By orthogonality, for all  $\ttt\in\TTT\backslash \SSS$, we have
\[
\norm{y+\sum_{m\in\NN} z_m(\ttt)}^2=1+\eta.
\]
and
\begin{equation}
  \label{eq:tail}
  \norm{\sum_{m=N}^\infty z_m(\ttt)}=\sqrt{\eta 2^{-N+1}}.
\end{equation}
Since
\[
\text{\textsl{weak}}\lim_{\ttt{\longrightarrow}\sss_n}  \left[y+\sum_{m\in\NN}
z_m(\ttt)\right]=
y+\sum_{m\not=n}z_m(\sss_n)
\]
we get again by orthogonality that the limit vector has  length
$\sqrt{1+\eta(1-2^{-n})}$. We are then ready
to define
\[
x(\ttt)=
\begin{cases}
\dfrac1{\sqrt{1+\eta}}
\left[y+\sum_{m\in\NN}z_m(\ttt)\right]&\ttt\in \TTT\backslash \SSS\\
\dfrac1{\sqrt{1+\eta(1-2^{-n})}}
\left[y+\sum_{m\not=n}z_m(\ttt))\right]&\ttt=\sss_n
\end{cases}
\]
and
\[
\mu(\sss_n)=\dfrac{\sqrt{1+\eta}}{\sqrt{1+\eta(1-2^{-n})}},
\]
where we note that $\mu(s_n)\nearrow 1$ as $n\longrightarrow\infty$. The condition (i) is now clearly met, and we have from \eqref{eq:tail} that
\[
\norm{y-x(\ttt)}\leq\left(1-\frac1{\sqrt{1+\eta}}\right)+\norm{\sum_{m\in\NN}z_m(\ttt)}=1-\frac1{\sqrt{1+\eta}}+\sqrt{\eta}<\delta
\]
for any $\ttt\in\TTT\backslash\SSS$. In fact, this inequality readily extends
to all of $\TTT$, proving that $x(\ttt)$ satisfies (ii). To prove that it
satisfies (iii), fix $\ttt_0\in\TTT\backslash\SSS$ and $\delta>0$. For an $N$
to be
determined later, choose a neighborhood $U_N$ of $\ttt_0$ such that 
$
  U_N\cap\SSS\subseteq\{s_n\mid n\geq N\}
$.
For $\ttt\in U_N\backslash S$ we get
\[
\norm{x(\ttt_0)-x(\ttt)}\leq
2\frac{\sqrt{\eta2^{-N+1}}}{\sqrt{1+\eta}}+
\frac1{\sqrt{1+\eta}}\sum_{m=1}^{N-1}{\norm{z(\ttt_0)-z(\ttt)}}
\]
by \eqref{eq:tail} again. A similar, slightly more complicated computation gives
that for $\ttt\in U_N\cap S$, 
\begin{equation}\label{eq:iiiesti}
\norm{x(\ttt_0)-x(\ttt)}\leq
2\frac{1-\mu(s_N)}{\sqrt{1+\eta}}+
2\frac{\sqrt{\eta2^{-N+1}}}{\sqrt{1+\eta}}+
\frac1{\sqrt{1+\eta}}\sum_{m=1}^{N-1}{\norm{z(\ttt_0)-z(\ttt)}},
\end{equation}
where we have used that the missing term in $x(\omega)$ in this case is not
among the first $N-1$ and that $s(\mu_n)\geq s(\mu_N)$ when $n\geq N$. Thus the estimate in \eqref{eq:iiiesti} holds throughout $U_N$. We can choose $N$ such
that the two first terms above sum to no more than $\epsilon/2$, and
then the norm continuity of each $z_1,\dots,z_{N-1}$ on $U_N$ to find a
neighborhood $\ttt\in V\subseteq U_N$ on which the last term is also 
bounded by $\epsilon/2$.

To prove (iv) it suffices to consider the case of sequences
$\ttt_k\longrightarrow\sss_{n_0}$. We may even reduce to the two cases of
$\{\ttt_k\}\cap \SSS=\emptyset$ and $\{\ttt_k\}\subseteq\SSS$. In the first
case, we note that $\ttt\mapsto \inner{z_m(\ttt)}{\myw}$ can be
extended with zero to a continuous function on $\TTT$ for each $m$, uniformly
bounded by $\sqrt{\eta 2^{-m}}\norm{\myw}$. On $\TTT\backslash\SSS$
we have
\[
\inner{x(\ttt)}{\myw}=\frac1{\sqrt{1+\eta}}\left(
\inner{y}{\myw}+\sum_{m\in\NN}\inner{z_m(\ttt)}{\myw}\right)
\]
and the expression on the right is continuous on all of $\TTT$ by the
above.
Hence
\[
\wlim_{k\longrightarrow\infty}x(\ttt_k)=
\frac1{\sqrt{1+\eta}}\left(
{y}+\sum_{m\in\NN, m\not=n_0}{z_m(\sss_{n_0})}\right)=
\mu(\sss_{n_0})x(\sss_{n_0}).
\]

In the second case, we write $\ttt_k=\sss_{n_k}$ and let $\epsilon>0$
be given. Because all values
$\sss_n$ are distinct, we can find, for some $N> n_0$ to be
determined, a neighborhood $U_N'$ of $\sss_{n_0}$ with the property
\begin{equation}
  \label{eq:noSmallNprime}
  U_N'\cap\SSS\subseteq\{s_{n_0}\}\cup\{s_n\mid n\geq N\}
\end{equation}
Choose $K$ with $\sss_{n_k}\in U_N'$ for all $k \geq K$. As in the proof of
(iii), we get that
\begin{eqnarray*}
\lefteqn{\inner{\mu(\sss_{n_0})x(\sss_{n_0})-x(\sss_{n_k})}{\myw}}&&\\
&\leq&
2\left(\frac1{\sqrt{1+\eta}}-\frac1{\sqrt{1+\eta(1-2^{-N})}}\right)\norm{\myw}+
2\frac{\sqrt{\eta2^{-N+1}}}{\sqrt{1+\eta}}\norm{\myw}\\
&+&\frac1{\sqrt{1+\eta}}
\sum_{m=1}^{n_0-1}\norm{{z_m(\sss_{n_0})-z_m(\sss_{n_k})}}\norm{\myw}+
\frac{\inner{z_{n_0}(\sss_{n_k})}{\myw}}{\sqrt{1+\eta}}\\
&+&\frac1{\sqrt{1+\eta}}{\sum_{m=n_0+1}^{N-1}\norm{z_m(\sss_{n_0})-z_m(\sss_{n_k})}}\norm{\myw}
\end{eqnarray*}
Choose first $N$ to make the two first terms less than $\epsilon/5$,
a corresponding $K$, and  then $K_1\geq K$ so that  the last three
terms each are
less than $\epsilon/5$ for $k\geq K_1$.
\end{demo}

We now consider $\A = C(\TTT,\KKK)$, where $\KKK$ is considered as  the compact
operators on our Hilbert space $H$. By standard identifications, 
$\A\dd\myz = \{f:\TTT\longrightarrow \BBB(\HHH)\mid f\text{ is bounded}\}$.
We
define  $r_x$  in  $\A\dd\myz$ 
by 
defining it at each  $\ttt$  in  $\TTT$ to be the projection on the span of
$x(\ttt)$. 

\begin{lemma}{zclrisr} With $x$ as in Lemma \ref{imprConstr}, for
  $\delta<1$, we have $\overline{r_x}\myz = r_x$.
\end{lemma}
\begin{demo}
Clearly $r_x=r_x\myz\leq \overline{r_x}\myz$, to reach a contradiction assume 
$r_x<\overline{r_x}\myz$. This implies that $(\overline{r_x}-r_x)\myz$
dominates a minimal projection, so some pure state of
$C(\TTT,\KKK)$ will evaluate to one on  $\overline{r_x}-r_x$. Such a
pure state will be given by
\[
f\mapsto \inner{f(\ttt_0)\myw_0}{\myw_0}
\]
for some choice of $\ttt_0\in \TTT$ and  unit vector $\myw_0\in \HHH$,
so we have, using the identification of $\A\dd\myz$ as operator-valued
functions on $\TTT$, that
\[
 \overline{r_x}(\ttt_0)\myw_0=\myw_0\qquad r_x(\ttt_0)\myw_0=0.
\]
By the definition of $r_x$, the second property implies orthogonality
of $\myw_0$ and $x(\ttt_0)$. Let $p,q\in \{f:\TTT\longrightarrow
\BBB(\HHH)\mid f\text{ is bounded}\}$ be given as projections onto
\[
\operatorname{span}\{\myw_0,y\}\qquad \operatorname{span}\{px(\ttt)\},
\]
respectively (recall that  $y$ is the distinguished  unit vector in $H$). Since $p$
is a  constant
projection, $p={p'}\myz$ where $p'\in\A$. Note also that $q$ has constant rank
one because
\[
\norm{px(\ttt)}\geq\norm{y}-\norm{p(y-x(\ttt))}\geq 1-\delta.
\]

In fact, we are going to prove that $q$ varies norm continuously with
 $\ttt$. Before we do so, let us show how this leads to the desired
 contradiction.
As above, $q=q'\myz$ where $q'\in\A$.
Since we have arranged that $p-q(\ttt)$
annihilates  $x(\ttt)$  for all  $\ttt$, we get  $(p-q)r_x=0$. This
entails that $r_x\leq 1-(p'-q')$, and since the larger projection here
is closed, even that
$\overline{r_x}\leq 1-(p'-q')$. Thus
$(p-q)\overline{r_x}=0$, leading to the contradiction 
\[
0=(p-q(\ttt_0))\overline{r_x}(\ttt_0)\myw_0=(p-q(\ttt_0))\myw_0= \myw_0.
\]

To see that $\ttt\mapsto q(\ttt)$ is norm continuous, note that by
what has already been said, 
\[ 
\myw(\ttt)=\frac1{\norm{p(x(\ttt))}}px(\ttt)
\]
is always defined. We have that $p$, being of finite rank, is
\wea-norm continuous, so we conclude from Lemma
\ref{imprConstr}(iii)-(iv) that
\[
\lim_{\ttt\longrightarrow\ttt_0}
px(\ttt)=\begin{cases}\mu(\ttt_0)px(\ttt_0)&\text{if }\ttt_0\in\SSS\\
px(\ttt_0)&\text{if }\ttt_0\not\in\SSS,\end{cases}
\]
showing that $\ttt\mapsto \myw(\ttt)$ is a norm continuous function in
$H$. When now $v$ is any unit vector, we have
\begin{eqnarray*}
\lefteqn{\norm{(q(\ttt) - q(\ttt_0))(v)}}\\
 &=& \norm{\inner{v}{ \myw(\ttt)} \myw(\ttt) - \inner{v}{ \myw(\ttt_0)}
 \myw(\ttt_0)}\\
& \leq&
\norm{\inner{v}{ \myw(\ttt)-\myw(\ttt_0)} \myw(\ttt)} +  \norm{\inner{v}{ \myw(\ttt_0)}
  (\myw(\ttt) - \myw(\ttt_0))}\\
&\leq&
\norm{\myw(\ttt)-\myw(\ttt_0)}^{1/2} \norm{v}^{1/2} \norm{\myw} +
\norm{v}^{1/2}\norm{\myw(\ttt_0)}^{1/2}\norm{\myw(\ttt)-\myw(\ttt_0)} \\
&\leq&
\norm{\myw(\ttt)-\myw(\ttt_0)}^{1/2} + \norm{\myw(\ttt)-\myw(\ttt_0)},
\end{eqnarray*}
and $q$ is norm continuous.
\end{demo}

\subsection{Dense nonregularity }\label{dn}
  We are going to prove that  $1 - \overline{r_x}$  is
dense and nonregular.  

\begin{propo}{dense}
  With $x$ as in Lemma \ref{imprConstr}, for any $\delta$, we have:
  \begin{rolist}
    \item $r_x$ is discontinuous in norm at any $\sss\in\SSS$.
\item $\myu-\overline{r_x}$ is dense.
  \end{rolist}
If furthermore $\delta<1/5$,we have
\begin{rolist}\addtocounter{rocount}{2}
\item $\myu-\overline{r_x}$ is not regular.
\end{rolist}
\end{propo}
\begin{demo}
For (i),  apply $r_x$ to  $x(\sss)$ and note that, as $\ttt$
approaches  $\sss$, 
\[
\ttt\mapsto r_x(\ttt)x(\sss)=\inner{x(\sss)}{x(\ttt)} {x(\ttt)}
\]
cannot converge in norm  to a unit vector
because $\mu(\sss)<1$.

For (ii), it suffices to show that  $\overline{r_x}$ could not dominate a
nonzero positive element  of  $\A$.    If it did, then by functional calculus, it would also
dominate an element $b$ with the property that   $b(\ttt)$ was a
projection on a nonempty open set  $U$  of  $\TTT$.
 By
Lemma \ref{zclrisr}, $r_x$ would dominate the canonical image of $b$
in $\myz\A\dd$, and since 
$r_x(\ttt)$ is rank 1 then
$b(\ttt)=r_x(\ttt)$  on 
$U$.  This contradicts (i) as  $b$  is norm continuous at any
$\sss_n\in U$.

For (iii), let $\sss\in\SSS$, and  define a constant element
$c\in\A$  by  $c(\ttt)=r_x(\sss)$  for all  $\ttt$.    Since
$\overline{r_x}\myz=r_x$ by Lemma \ref{zclrisr},
\[
        (1-\overline{r_x}(\ttt))c(\ttt)\myz=(1-r_x(\ttt))r_x(\sss)=(r_x(\sss)
        - r_x(\ttt))r_x(\sss),
\]
so for any $\ttt\in\TTT$ and any unit vector $v\in H$,
\begin{eqnarray*}
\lefteqn{\norm{(r_x(\sss) - r_x(\ttt))v}=}&&\\
&&\norm{\inner{v}{x(\sss)}x(\sss) - \inner{v}{x(\ttt)}x(\ttt)}\leq
2\delta+\norm{\inner{v}{y}y - \inner{v}{y}y}+2\delta\leq \frac45.\\
\end{eqnarray*}
This gives, using \cite[4.3.15]{gkp:cag},
\begin{eqnarray*}
\norm{(1-\overline{r_x})c}&=& \norm{(1-\overline{r_x})c\myz}\\
&=&
\sup_{\ttt,\norm{\myw}\leq 1}\norm{(r_x(\sss) - r_x(\ttt))r_x(\sss)\myw}\\
&\leq&
\sup_{\ttt,\norm{v}\leq 1}\norm{(r_x(\sss) - r_x(\ttt))v}\leq\frac45\\
&<& 1=\norm{c}=\norm{\overline{\myu-\overline{r_x}}c}, 
\end{eqnarray*}
and
 $1-\overline{r_x}$  is not regular,
by Theorem \ref{regumain}(i).
\end{demo}

\begin{theor}{denseopen}
Let $\cantor$ be the Cantor set. The \cstar-algebras $C([0,1],\KKK)$ and
$C(\cantor,\KKK)$ both have open, dense and nonregular projections. The latter
\cstar-algebra is $AF$.
\end{theor}
\begin{remar}{reguof}
  The reader may wonder whether $\myu-r_x$ is regular. In fact, since
  it dominates $\myu-\myz$ by construction, and since every ideal of
  $C(\TTT,\KKK)$ has diffuse states, our Proposition \ref{coz} below
  will show that it is indeed regular.
\end{remar}

\subsection{Dense nonregularity in simple \cstar-algebras}\label{dns}

We now proceed to show the existence of a dense open nonregular
projection in $\MMM\otimes\KKK$. Our proof is based on pushing the
projection constructed in $C(\cantor)\otimes\KKK$, with $\cantor$ the
Cantor set, in the previous
section into $\MMM\otimes\KKK$. This will prove to be possible
because, as noted in \cite{ja}, $C(\cantor)$ sits inside $\MMM$ as a
diagonal MASA on which  the ``diagonal compression'' map
$E:\MMM\longrightarrow C(\cantor)$ is a faithful conditional
expectation (a projection of norm one).

One should note right away that since, as we shall see in Remark \ref{hownotto}, regularity is not a hereditary property, one can not in general
conclude the existence of a nonregular projection in a \cstar-algebra 
$\A$ from the existence of such a projection in a subalgebra of $\A$. 
The proposition below, a collection of well-known results, compiles
the information we shall need to make such an argument work.

\begin{propo}{Tomi}
There is a faithful conditional expectation $E$ from the
$2^\infty$ UHF algebra
$\MMM$  onto its diagonal MASA  $C(\cantor)$. Then $E$ induces a map
\[
E\otimes \id:\MMM\otimes\KKK \longrightarrow C(\cantor)\otimes\KKK
\]
which is again a faithful conditional expectation.
Furthermore, the normal extension
\[
(E\otimes\id)\dd : (\MMM\otimes\KKK)\dd \longrightarrow (C(\cantor)\otimes\KKK)\dd 
\]
is also a conditional expectation.
\end{propo}
\begin{demo}
The existence of $E$, which is faithful as it is trace-preserving, is
noted in \cite{ja}.
The tensor product then exists as $E$ is a completely positive map,
and by \cite{jt:pnow}, is will be a faithful conditional expectation.
That also $(E\otimes \id)\dd$ is a conditional expectation follows by 
routine duality arguments.
\end{demo}

Denote by  $\iota$ the inclusion map of  $C(\cantor)$  into
$\MMM$.  Since  $\iota$  is a $*$-monomorphism, so is the canonically
induced map
\[
(\iota\otimes\id)\dd:(C(\cantor)\otimes\KKK)\dd \longrightarrow
(\MMM\otimes\KKK)\dd.
\]

\begin{propo}{odnrinsim}
Let $p$ be an open dense nonregular projection for
$C(\cantor)\otimes\KKK$ . Then the projection $(\iota\otimes\id)\dd(p)$  is an open, dense, nonregular projection for  $\MMM\otimes\KKK$.
\end{propo}
\begin{demo}
Since $p$ is open, dense and nonregular, there is a nonzero element  $d_0  \in  C(\cantor)\otimes\KKK$  such that
$\norm{d_0} = 1 > \norm{pd_0}$.
Since $(\iota\otimes\id)\dd$ is a $*$-monomorphism,  $(\iota\otimes\id)\dd(p)$  is an open projection.  Suppose
we can prove that it is dense for  $\MMM\otimes\KKK$.  Then it must be nonregular since
$\norm{d_0} = 1 > \norm{pd_0}$  and hence
$1 = \norm{(\iota\otimes\id)(d_0)} > \norm{(\iota\otimes\id)\dd (p) (\iota\otimes\id)(d_0)} $.

        It remains to show that  $(\iota\otimes\id)\dd(p)$  is dense  for  $\MMM\otimes\KKK$.  Let
$b$  denote a positive element of  $\MMM\otimes\KKK$  such that
$(\iota\otimes\id)\dd(p)b=0$.  Using the properties of conditional expectations,
\[
0 = (E\otimes\id)\dd((\iota\otimes\id)\dd(p)b) = (\iota\otimes\id)\dd(p)(E\otimes\id)\dd(b) =
(\iota\otimes\id)\dd(p(E\otimes\id)(b)).
\]
But  $(\iota\otimes\id)\dd$ is injective, so this means that  $p (E\otimes\id)(b) = 0$.
Since $p$  is dense for $C(\cantor)\otimes\KKK$  and  $E\otimes\id(b)$
is in $C(\cantor)\otimes\KKK$, $E\otimes\id(b) = 0$.  By Proposition \ref{Tomi} above,  $E\otimes\id$  is faithful, so  $b \geq 0$ and $E\otimes\id(b) = 0$  imply that  $b =
0$.  Thus  $(\iota\otimes\id)\dd(p)$  is dense, and we are done.
\end{demo}

\begin{theor}{denseopensimple}
  There is a dense open projection for the stabilized $2^\infty$ UHF
  algebra $\MMM\otimes\KKK$.
\end{theor}
\begin{demo}
  Combine Proposition \ref{odnrinsim} and  Theorem \ref{denseopen}.
\end{demo}
\begin{remar}{hownotto}
  Consider
\[
\A=\left\{(A_n)_{n\in\ZZ}\in(\matrM_2)^\ZZ\left|
\lim_{n\longrightarrow\pm\infty}A_n\text{ exists}
\right.\right\}
\]
and its subalgebra
\[
{\mathfrak B}=\left\{(A_n)_{n\in\ZZ}\in\A\left| \lim_{n\longrightarrow+\infty}A_n=
\lim_{n\longrightarrow-\infty}A_n\right\}\right.
\]
As the  irreducible representations of ${\mathfrak B}$ (all 2-dimensional) can be
naturally parametrized over   $\ZZ\cup\{\infty\}$, we can specify an
open projection $p$ in ${\mathfrak B}\dd$ by
\[
p\myz(n)=\begin{cases}\smalltwobytwo 1000&n\in \NN_0\\
\smalltwobytwo {1/2}{1/2}{1/2}{1/2}&n\in -\NN\\
\smalltwobytwo 0000&n=\infty
\end{cases}
\]
Since one finds that $\overline{p}\myz$ differs only from $p\myz$ by
having $\overline{p}\myz(\infty)=\smalltwobytwo1001$, one can get by
considering the element $b\in{\mathfrak B}$ which is constantly
$\smalltwobytwo0001$ that $p$ is not regular. 
However, $p$ considered as an element of $\A$ has a closure which lies 
in ${\mathfrak B}$ itself, and thus must be regular. %
\end{remar}

\subsection{Nonregularity in $\MMM$}\label{ninM}
At present, we do not see a way to amend our construction to prove
that there is a dense open nonregular projection for the $2^\infty$
UHF algebra $\MMM$ itself, but only for its stabilized version. The
difference in this respect between $C(\TTT,\matrM_n)$ and
$C(\TTT,\KKK)$ may seem to indicate that this is due to some deeper
phenomenon.

However, we can supply an example of a nonregular open projection
associated to $\MMM$, and we supply the details here as we feel the
methods may be of independent value.

\begin{lemma}{projsczP}
If $p\in\A\dd$ is regular and dominated by $\myz$, then
\[
\overline{\PPA{p}}=\PPA{\pbar},
\]
where the closure on the left hand side is relative to $P(\A)$.
\end{lemma}
\begin{demo}
We first apply the assumption that $p\leq\myz$ to prove that
$\cco{P(p)\cup\{0\}}=F(p)\cl$. Since
$F(p)\cl$ is a \weasta\
closed convex set containing $\PPP{p}$, the inclusion from left to right is
clear. By the double polar theorem (as it is found in, e.g.,
\cite[IV.1, proposition 3]{nb:evt})
applied to the real 
spaces $(\A\sa,\norm{\cdot})$ and $(\A^*\sa,\weasta)$ in
duality
(cf.\ \gkp{3.1.1}) and the fact that $0\in
F(p)$, the other inclusion is equivalent to
$\PPP{p}^{\circ\circ}\supseteq F(p)^{\circ\circ}$, which follows
from
\begin{eqnarray*}
P(p)^\circ&=&\{a\in\A\sa|\forall \phi\in P(p):\phi(a)\geq -1\}\\
&=&\{a\in\A\sa|\forall \phi\in P(\A):\phi\dd(pap)\geq -1\}\\
&=&\{a\in\A\sa|\forall \pi\in \hat{\A}:\pi\dd(pap)\geq -\myunit\}\\
&=&\{a\in\A\sa|pap\myz\geq -\myz\}\\
&=&\{a\in\A\sa|pap\geq -\myunit\}\\
&\subseteq &F(p)^\circ.
\end{eqnarray*}
For the inclusion from left to right in the second equality, let
$a\in\A\sa$ with $\phi(a)\geq-1$ for all $\phi\in P(p)$ be
given and set
$\alpha=\psi(p)$ for a given $\psi\in P(\A)$. If $\alpha=0$, clearly
$\psi\dd(pap)\geq-1$. If $\alpha>0$, one uses that
$\psi'=\alpha^{-1}\psi(p\cdot p)$ is an element of $P(p)$.
The assumption on $p$ is used in the inclusion from left to right of
the last equation; if $pap\myz\geq-\myz$, we have
\[
pap=pap\myz\geq -\myz\geq-\myunit.
\]
 
When $\PPP{p}\cl$ denotes the closure of $\PPP{p}$ relative to $P(\A)$,
it is clear that $\PPP{p}\cl\subseteq\PPP{\pbar}$. Using
the fact that $p$ is regular along with the above observation, we get
\[
\cco{\PPP{p}\cup\{0\}}=F(p)\cl=F(\pbar).
\]
By \cite[Appendice B 14]{jd:clr}, the extremal points of
$F(\pbar)$ are contained in $(\PPP{p}\cup\{0\})\cl$, and as
$F(\pbar)$ is a face, these are exactly
$\PPP{\pbar}\cup\{0\}$. We conclude that $\PPP{\pbar}\subseteq
\PPP{p}\cl$, as required.
\end{demo}

\begin{lemma}{qreg_iff_qzreg}
Assume that  $q$  is an open projection of $\A$. Then $q$ is regular if and
only iff $q\myz$ is. In this case, $\overline{q\myz}=\overline{q}$.
\end{lemma}
\begin{demo}
Assume that $q$ is regular.  Since $q$ is open,   $q$  is a limit of a net of
compact projections  ${q_k}$ dominated by
$q$.
  Since  $F(\myz q_k)\cl= F(q_k)$, then $F(q)$ is contained in the closure of
  the
union of the $F(\myz q_k)$.  Thus  $F(q)$  is contained in $F(\myz q)\cl$, and
\[
F(\overline{q\myz})\subseteq
F(\overline{q})=\overline{F(q)}=\overline{F(q\myz)} \subseteq
F(\overline{q\myz}).
\]

In the other direction, suppose that  $\myz q$  is regular, and  let  $r$
be the closure of  $\myz q$.  Then  $\myz r$  dominates  $\myz q$.  Since  $r$
 is closed and
$q$  is open, \cite[4.3.15]{gkp:cag} shows that  $r$  dominates  $q$.  Thus
the closure of
$q$  is dominated by  $r$, hence it equals  $r$. We conclude that
$\overline{q\myz}=\overline{q}$. But  $\myz q$ regular means that for any  a
in  $\A$,  $\norm{\myz qa}=\norm{ra}$,
but also
$\norm{\myz qa}\leq\norm{qa}\leq\norm{ra}$  because of the ordering of  $\myz
q$, $q$ and $r$.
Thus
$\norm{qa}=\norm{ra}$,  so  $q$  is regular.
\end{demo}

Now for an open projection $q$, we can look at $q\myz$ and use Lemma
\ref{qreg_iff_qzreg} to check for regularity.

\begin{propo}{clopencriterion}
If $q$ is open and regular, and if $\overline{P(\myz q)}=P(\myz q)$, then
$q\in
M(\A)$.
\end{propo}
\begin{demo}
Since   $q$ is regular,  $\myz q$  is regular by Lemma
\ref{qreg_iff_qzreg} above.  By
Lemma \ref{projsczP}, regularity of $\myz q$ means that $P(\myz
q)\cl=P(\overline{q})$, implying that $q\myz =\overline{q}\myz$.  As in
\cite[4.3.15]{gkp:cag}, because
$q$  and  $q\cl$  are both semicontinuous, we conclude that $q=\overline{q}$.
Since  $q$  is
already open, it lies in  $M(\A)$  by \cite[2.2]{caagkpjt:mc}.
\end{demo}

We now apply  the results above to the ``perfection'' $\A_c$  of
a \cstar-algebra $\A$ introduced in \cite{caafws:pc}. This \cstar-algebra is
defined  as
\[
\{a\in \A\dd\myz\mid a,a^*,aa^* \text{ are \textsl{weak}*-continuous on }
P(\A)\cup\{0\}\},
\]
and has subsequently been studied in \cite{caajagkp:dspc} and
\cite{cjkb:scsw}. This notion ties in with regularity as follows:

\begin{corol}{perfectionandregularity}
  If $p$ is an open and regular projection of $\A$, and if $p\myz\in \A_c$,
  then
  in fact $p\in M(\A)$.
\end{corol}
\begin{demo}
By definition of  $\A_c$,  elements of  $\A_c$  are continuous on the pure
states of  $\A$. Hence  $P(p\myz)\cl=P(p\myz)$, and Proposition
\ref{clopencriterion} applies.
\end{demo}

We are now ready to prove the existence of a nonregular projection in $\MMM$. Before we give the proof, let us  review the
following key notions from \cite{caajagkp:esc} and \cite{caajagkp:dspc}:

\begin{defin}{excision}
A net $x_\lambda$ of elements in $\A^+$ with $\norm{x_\lambda}=1$
\emph{excises} the state $\phi\in S(\A)$ if
\[
\norm{x_\lambda(a-\phi(a))x_\lambda}\longrightarrow 0, \qquad\forall a\in \A.
\]
\end{defin}

\begin{remar}{checklessinexcision}
To check that a net  $x_\lambda$ in $A^+$ with $\norm{x_\lambda}=1$ 
excises a state $\phi$ it  suffices to check the convergence for each
$a$  in a dense subset of  $\A$. 
\end{remar}
\begin{defin}{diffusedefs}
Let $\A$ be a \cstar-algebra. 
\begin{rolist}
  \item (\cite[2.1]{caajagkp:dspc}) A sequence $a_n$ of $\A$ is \emph{diffuse}
    if for every net
    $\phi_\lambda$ in $P(\A)\cup\{0\}$, converging \weasta\ in
    $P(\A)\cup\{0\}$, we have
\[
\lim_{\lambda,n}\phi_\lambda(a^*_na_n+a_na_n^*)=0.
\]
\item (\cite[3.1]{caajagkp:dspc}) An orthogonal, positive, norm one
  sequence $a_n$ in $\A$ is \emph{truly diffuse} if for any increasing
  sequence $n_k$ in $\NN$, the sequence 
\[
\sum_{j=n_k}^{n_{k+1}-1}a_j
\]
is diffuse.
\item (\cite{caajagkp:esc}) A quasi-state $\phi\in Q(\A)$ is
  \emph{diffuse} when $\phi(\myz)=0$.
\end{rolist}
\end{defin}

\begin{theor}{nonregularinUHF}
The $2^\infty$ UHF algebra $\MMM$ contains a nonregular open projection.
\end{theor}
\begin{demo}
        By \cite[3.7]{caajagkp:dspc} there is a factor state  $\phi$  of type
        II$_\infty$, a
sequence $\{p_n\}$ of orthogonal projections in  ${\MMM}$  and a dense
        sequence  $\{a_n\}$ 
in  ${\MMM}$  such that for $1 \leq k \leq n$ and  $m>n$,
\begin{rolist}
  \item $p_n(a_k - \phi(a_k))p_n = 0$
\item $p_ma_kp_n = 0 = p_na_kp_m.$
\end{rolist}

We first show that $\{p_n\}$  is a truly diffuse sequence.  Let  $\{n_k\}$  be
an
increasing sequence of natural numbers and define  $c_k =
\sum_{j=n_k}^{n_{k+1}-1}p_j$;  by Definition \ref{diffusedefs}(ii) we need to
show that  $\{c_k\}$  is a
diffuse sequence.  Using 
\cite[2.14]{caajagkp:dspc} we can conclude that  $\{c_k\}$  is a diffuse
sequence if it excises
a diffuse state.  Since  $\phi$  is of type II$_\infty$, it is diffuse.  Since
$\norm{c_k}= 1$  for all $k$,  by Remark \ref{checklessinexcision} we only
need to check the
excising condition on the dense set  $\{a_k\}$, and here it is
immediately verified for each  $a_k$  by (i) and (ii) above.  Thus  $\{p_n\}$
is
a truly diffuse sequence.

Set  $p = \sum_{n=1}^\infty{p_n}$, where the sum is taken in
${\MMM}\dd$.
If  $p$  were in  ${\MMM}$,
then by Dini's theorem the sequence of partial sums
$\{\sum_{n=1}^k{p_n}\}$,   would have to converge to  $p$  in norm, and
that is impossible. 
Thus  $p$  does not lie in  ${\MMM}$.

        We next show that $p\myz$ lies in ${\MMM}_c$.  Since
$p\myz$ is a projection and ${\MMM}$ is separable, it sufffices
to assume that $\{\psi_n\}$ is a sequence of pure states of
${\MMM}$ that converges to a pure state $\psi$ of $\MMM$
and to show that $(\psi_n - \psi)(p)\longrightarrow 0$.  Let $\epsilon
> 0$ be given. Choose $n_0$ such that $\psi(\sum_{j=n_0}^\infty{p_j})<
\epsilon/3$.  This is possible since $\sum_{n=1}^\infty{p_n}$ is
\weasta\ convergent in ${\MMM}\dd$ and $\psi$ is \weasta\
continuous on $A\dd$.  Since $\{p_n\}$ is a truly diffuse sequence,
\cite[3.2]{caajagkp:dspc} allows us to find $n_1 > n_0$ such that for
$j\geq n_1$, $\psi_j(\sum_{n=n_1}^\infty{p_n})<\epsilon/3$.  Choose $n_2
> n_1$ such that for $j > n_2$, $|(\psi_j -
\psi)(\sum_{i=1}^{n_1-1}{p_i}) < \epsilon/3$.  Now for $j > n_2$,
\[
{\left|\left(\psi_j - \psi\right)\left(\sum_{i=1}^\infty{p_i}\right)\right|<}
\left|\left(\psi_j - \psi\right)\left(\sum_{i=1}^{n_1-1}{p_i}\right)\right|
+\left|\psi_j\left(\sum_{n=n_1}^\infty{p_n}\right)\right| +
\left|\psi\left(\sum_{n=n_1}^\infty{p_n} \right)\right| <
\epsilon.
\]
This shows that  $p\myz$  is in  ${\MMM}_c$. Since  $p$  was not in
$M({\MMM})={\MMM}$, then  $p$  is not regular by  Corollary
\ref{perfectionandregularity}.
\end{demo}

\section{Automatic regularity of large projections}
In this section, we present a few positive results on regularity that
we found while trying to settle the general questions described in the
introduction. It is our hope that they can be used to work around some of
the complications that the existence of dense open projections lead
to.

%
%
%
%
%
%
%
%
%
%
%
%
%
%
%
%
%
%
%
%
%
%
%
%
%
%
%
%
%
%
%
%
%
%
%
%
%
%
%
%
%
%
%
%
%
%
%
%
%
%
%
%
%
%
%
%
%
%
%
%
%
%
%
%
%
%
%
%
%
%
%
%
%
%
%
%
%
%
%
%
%
%
%
%
%
%
%

\subsection{``Bottom  up'' regularity}

The following lemma follows directly from the definitions and the fact
that central projections are regular. Nevertheless, it
plays a key role in establishing our more  profound results at the end
of the section.

\begin{lemma}{char_reg_new}
Let $p$ be a projection of the \cstar-algebra $\A$. 
\begin{rolist}
\item $p$ is regular and dense if and only if
\[
\forall a\in \A: \norm{ap}=\norm{a}
\]
\item If $p$ dominates a regular and dense projection, then $p$ is
  also regular and dense.
\item If $p$ dominates a central and dense projection, then $p$ is
 regular and dense.
 \end{rolist}
\end{lemma}

As an example of the strength of this form of reasoning, note that it
gives a short alternative proof of \cite[3.4]{cplz:opccr} since when
$\KKK$ is an essential ideal in $\A$ and $p$ is dense and open in
$\A\dd$, then $p$ must dominate the cover of $\KKK$, which is 
dense and central.

Before we move on to further consequences, we need a few preliminaries:

\begin{remar}{onideals}
If  $\I$  is a closed ideal of  $\A$, then there is a central open
projection  $x$ in  $\A\dd$
such that  $\I = \A \cap x\A\dd$. In this setting,  and $\I\dd$ is isometrically isomorphic to
$x\A\dd$ (see \cite[3.10.7]{gkp:cag}). Further,  $(\A/\I)\dd$ is isometrically
isomorphic to $(1-x)\A\dd$.
The first isomorphism respects $\myz$ in the sense that
\[
x\myz_\A = \myz_\I
\]
when $\myz_\I$ denotes the sum of the minimal projections in
$\I\dd$ and $\myz_\A$ denotes  the sum of the minimal projections in $\A\dd$
(see \cite[3.13.6(iii)]{gkp:cag}). 
\end{remar}

The notion of a \emph{scattered} \cstar-algebra from \cite{hejI} will also be
useful. Here, a \cstar-algebra is defined to be scattered if 
no state of $\A$ is diffuse, cf.\ Definition \ref{diffusedefs}(iii). 
 By \cite[2.2]{hejI}  $\A$  is scattered precisely if 
$\myz = 1$.  We use these facts in the next proof.

\begin{propo}{coz}
  Let $\A$ be a \cstar-algebra  that has no nonzero scattered ideal.  Then any projection dominating $\myu-\myz$ is a dense regular projection.
\end{propo}
\begin{demo}
By Lemma \ref{char_reg_new}(iii) it suffices to prove  that $\myu-\myz$
is dense.  Let  $x= 1 - (1-\myz)\cl$. If $x=0$, then $\myu-\myz$ is dense, so, to reach a contradiction, 
assume that $x\not=0$.  Since  $x$  is an open central
projection, $\I=\A\cap x\A\dd$ is a nonzero ideal. By hypothesis $\I$ is not
scattered, so after identification as explained above,  $\myz_\I < x$.  Thus,
by Remark \ref{onideals} again,
\[
0< (x - \myz_\I) \leq 1 - \myz_\A,
\] contradicting the definition of  $x$.
\end{demo}

\begin{corol}{scatanti}
    If  $\A$  is antiliminary, then any projection dominating $\myu - \myz$  is
dense and regular.
\end{corol}
\begin{demo}
  By \cite[3.2]{hejI}, any scattered \cstar-algebra is type I.  Since  $\A$  is
antiliminary, it has no nonzero type I ideals, so the conclusion follows
by Proposition \ref{coz}.
\end{demo}

\begin{propo}{411}  {Assume that  $\A$  is a \cstar-algebra with a faithful tracial
state  $\tau$.  If  $p$  is an  open projection in  $\A\dd$ such that  $\tau(p) = 1$,
then  $p$  is regular and dense.}
\end{propo}
\begin{demo} { Let  $x$  be the support projection of  $\tau$  in  $\A\dd$.  Since  $\tau$ is
a trace, it is unitarily invariant, and that immediately implies that  $x$
is a central projection.  Since  $\tau(p) = 1$, $p$ must dominate  $x$  by the
definition of  $x$  as the support projection of  $\tau$.  Since  $\tau$  is
faithful, $x$ is dense.  Since  $x$  is central, it is regular, so Lemma
\ref{char_reg_new}(ii) implies that  $p$  is regular and dense.}
\end{demo}

\begin{remar}{bottoprem}  {The situation in Proposition \ref{411} can arise in many ways.  For
example, let  $\MMM$  be the $2^\infty$ UHF algebra with trace  $\tau$.  Recursively
choose an orthogonal family  $\{p_n\}$  of projections in  $\MMM$  such that
$\tau(p_n) = 2^{-n}$.  (There are uncountably many distinct ways to do this.)
Then  $\sum\{p_n\}$, taken in  $\MMM\dd$, is an open projection with trace 1.  Even in
an algebra with no nontrivial projections, e.g. the reduced \cstar-algebra of
the free group on two generators, this same recursive construction can
take place, except that the projections  $p_n$  will be open in $\MMM\dd$, not
lying in  $\MMM$  itself.  We view this as a "bottom up" method of constructing
dense, regular open projections.  By contrast, the "top down" method of
Corollary \ref{cofinite_regular_ii} below shows that certain projections that are constructed by
subtracting closed projections from the identity are also regular and dense.}
\end{remar}

We end this section with the following useful lemma, in which we 
 tacitly invoke the isomorphisms from Remark \ref{onideals}.
We point out that one direction of (ii) below was already noted in \cite[3.5]{cplz:opccr}.

\begin{lemma}{reduclem}
 Let  $\A$  be a \cstar-algebra with an ideal  $\I$  whose central
cover in  $\A\dd$ is  $x$.
\begin{rolist}
\item A projection  $q$  in  $\A\dd$  is regular
and dense for  $\A$ 
if and only if  $xq$  is regular and dense for  $\I$  and  $(1-x)q$  is
regular and dense for
$\A/\I$.  
\item
{  If  $\I$  is essential, then  $q$  is regular and dense for  $\A$
if and only if  $xq$  is regular and dense for  $\I$.}
\end{rolist}
\end{lemma}
\begin{demo}\hide{
  The forward direction is trivial from Lemma \ref{char_reg_new}(i).  Now assume
  that  $xq$  is
regular and dense for  $\I$  and  $(1-x)q$  is regular and dense for
$\A/\I$.  Let
$b$  be a  norm $1$ element of  $\A$.  Suppose that  $\norm{bq} <
1$.  Then
$\norm{b}=1$ implies that either $\norm{xb}=1$ or  $\norm{(1-x)b}=1$.  If
$\norm{(1-x)b}=1$, then  $\norm{(1-x)bq}=1$  by the regularity and density of
$(1-x)q$, contradicting  $\norm{bq} < 1$.  Thus we can assume that
$\norm{xb}=1$.
Choose a positive norm one element  $a$  of  $\I$  such that  $\norm{ab} >
\norm{bq}$.
This is possible since  $x$  is the \weasta\ limit of 
elements of  $\I$ of norm less than one, hence  $xb$  is the \weasta\
 limit of elements of the form  $ab$. 
Since norm closed balls are also \weasta\ closed (by the definition of the
dual space norm), the set of all elements of the form  $ab$  can't all be
contained in a ball about  $0$  with radius strictly less than  $\norm{xb}=1$.
However,  $ab$  lies in  $\I$, so by regularity and density of  $xq$,
\[
\norm{bq}=
\norm{a} \norm{bq} \geq\norm{abxq}=\norm{ab} > \norm{bq},
\]
a contradiction, proving (i)  
}

{ By Lemma \ref{char_reg_new}, it suffices to show that if  $xq$  is regular and
dense for  $\I$, then  $q$  is regular and dense for  $\A$.  Let  $b$  lie in  $\A$
with  $\norm{b} = 1$.   It suffices to show that  $\norm{qb}$ = 1.
Since  $\I$  is essential, the map $A \mapsto xA$  has no kernel, hence is isometric.
Thus $\norm{xb} = 1$, and   the proof of (ii) now proceeds as the proof of (i).
}

\end{demo}

\subsection{``Top down'' regularity}

We end the paper by a collection of results which generalize
\cite[3.6]{cplz:opccr}. Many ingredients in the proof below were
indeed borrowed from that source. 

\begin{theor}{new43}  Let  $\{p_n\}_{n\in\mathcal N}$  be a countable
  family  of minimal projections in  $\A\dd$, with the property that
 with
\[
\mathcal F=\{n\in\mathcal N\mid p_n\in\A\},
\]
then  $\mathcal F$ is 
a finite  (possibly void) set.  Set  
\[
p=\bigvee_{n\in\mathcal F}p_n
\qquad
q=1-\bigvee_{n\in\mathcal N}p_n.\]
 Then $q$  is regular and $q\cl=1-p$.
\end{theor}
\begin{demo}
The proof will be in steps and will keep the notation above.  Other
symbols may be reused from one step to the next.

\noindent\textsc{Step 1:  Reduction to   $\mathcal{F}=\emptyset$.}
        Since each  $p_n$  for  $n$  in  $F$  is both open and closed, so is
        the
supremum by  \cite[2.5 \& 2.9]{caa:gswp}, hence $1-p$
is both open and closed, and consequently  a multiplier of  $\A$
(\cite[2.5]{gkp:awsct}).  Thus  $q\cl$ can be no larger
than $1-p$.  Both conclusions of the theorem will therefore follow if we
can show that for any positive, norm 1 element  $b$  of  $(1-p)\A(1-p)$,
$\norm{pb}=1$.  We may thus, without loss of generality, assume that  $p=0$,
i.e. that  $\mathcal F$  is void.

{\noindent\textsc{Step 2.  Reduction to cases.}
By \cite[6.2.7]{gkp:cag} there is a largest type I ideal  $\I$  of  $\A$  and
$\A/\I$  is antiliminary.  Let  $x$  be the central cover of  $\I$  in  $\A\dd$.  By
Lemma \ref{reduclem}(i) it suffices to show that  $xq$  is regular and
dense for  $\I$  and 
$(1 - x)q$  is regular and dense for  $\A/\I$.  We shall show in the next
paragraph that both  $xq$  and  $(1-x)q$ can be expressed as required in the
hypothesis  of the present 
theorem.  Thereafter it will suffice to demonstrate the theorem separately in the 
cases below.}

{Note that since  $x$  is central, a rank 1 projection in
$\A\dd$ lies under  $x$  or under  $1-x$.   Thus the projections  $\{p_n \mid n \in \NN\}$ are
partitioned into two subsets, those lying under  $x$  and those lying under
$1-x$.  Obviously a projection that lies in  $\I$  must lie in  $\A$.  Therefore
$xq$  is the complement of the supremum of a countable family of rank 1
projections in  $\I\dd$, none of which lies in  $\I$.  As for  $(1-x)q$, clearly it
is the complement of the supremum of a countable family of rank 1 projections
in  $(\A/\I)\dd$.  Since  $\A/\I$  is antiliminary, it can't contain any rank one
projections by \cite[6.1.7]{gkp:cag}.  Therefore both  $xq$  and  $(1-x)q$ can
be expressed as hypothesized.}

\noindent\textsc{Step 3.  The antiliminary case.}  Assume that  $\A$  is
antiliminary.
Since $q$ dominates $\myu-\myz$, we are done by Corollary \ref{scatanti}.

\noindent\textsc{Step 4.  The type I case}.  Assume that  $\A$  is type I.
        By \cite[6.2.11]{gkp:cag}, $\A$ contains an essential ideal  $\J$  that
        has
continuous trace. Arguing as above with Lemma \ref{reduclem}(ii) we can
pass to the case below.

\noindent\textsc{Step 5.  The continuous trace case}.  Assume that  $\A$  has 
continuous trace, and recall that the spectrum of  $\A$  is a locally compact Hausdorff space by
\cite[6.1.11]{gkp:cag}.  Since  $q$  dominates the
complement of the supremum $r$ of the central covers of the $\{p_n\}$, it
suffices to prove  that $r$ is a dense central projection. This
follows by a category argument, since $r$ is represented in the
spectrum of $\A$ as the complement of a countable set (namely the
central covers of the projections  $\{p_n\}$), which is still a dense set.
Of course countability of $\{p_n\}$  is crucial here.
\end{demo}

\begin{corol}{cofinite_regular_ii}
 Let  $\A$  be any \cstar-algebra, and suppose  $p\in\A\dd$  has finite
 codimension.  Then  $p$  is regular.
\end{corol}

\begin{corol}{NGCR}
 If  $\A$  contains no minimal projections, and $\{\phi_n\}_{n\in \NN}$
are pure states of  $\A$ with  $\{p_n\}_{n\in\NN}$  their corresponding support
projections
in  $\A\dd$, then 
\[
1-\bigvee_{n\in\NN}{p_n}
\]
  is
regular
and dense.
\end{corol}

\begin{remar}{scattered}
Note that when every ideal of $\A$ has a diffuse state, the
countability condition in the corollary above is unnecessary as any
projection of the form
\[
\myu-\bigvee_{i\in I}p_i,
\]
with $p_i$ minimal, dominates $\myu-\myz$ which is regular and dense
by Proposition \ref{coz}.

Indeed, it would be possible to 
strengthen the results above  further by combining these ideas.
We
refrain from this for the moment.
\end{remar}

\section*{Acknowledgment}
S{\o}ren Eilers was partially
    supported by the Carlsberg  Foundation in the initial stages of the work
    leading to the present paper.

\providecommand{\bysame}{\leavevmode\hbox to3em{\hrulefill}\thinspace}

\begin{center}
\parbox{6.5cm}
{\begin{center}{\sc
Department of Mathematics\\
University of California\\
Santa Barbara\\
California 93106\\
U.S.A.\\}
{\tt akemann@math.ucsb.edu}
\end{center}}
\makebox[0.5cm]{}
\parbox{6.5cm}
{\begin{center}{\sc
Matematisk Institut\\
K\o benhavns Universitet\\
Universitetsparken 5\\
DK-2100 Copenhagen \O\\
Denmark\\}
{\tt eilers@math.ku.dk}
\end{center}}\\
\end{center}
\end{document}